\documentstyle[12pt,a4,amsfonts]{article}

\newtheorem{Theorem}{Theorem}[section]
\newtheorem{Proposition}[Theorem]{Proposition}
\newtheorem{Lemma}[Theorem]{Lemma}

\newtheorem{Remark}[Theorem]{Remark}

\newtheorem{Definition}[Theorem]{Definition}

\def\RMN#1{\uppercase\expandafter{\romannumeral#1}}

\newcommand{\qed}{{\unskip\nobreak\hfil\penalty50\quad\null\nobreak\hfil{\bf
q.e.d.}\parfillskip0pt\finalhyphendemerits0\par\medskip}}
\newcommand{\proof}{\noindent{\it Proof.} \ }

\newcommand{\rank}{\mathop{\rm rank}\nolimits}

\newcommand{\proj}{\mathop{\rm Proj}\nolimits}

\renewcommand{\ker}{\mathop{\rm Ker}\nolimits}

\begin{document}

\title{
On finite generation of symbolic Rees rings of space monomial curves 
and existence of negative curves
}
\author{Kazuhiko Kurano and Naoyuki Matsuoka}
\date{
Dedicated to Professor Paul C. Roberts \\
on the occasion of his 60th birthday
}
\maketitle

\abstract{In this paper, we shall study finite generation
of symbolic Rees rings of the defining ideal of
the space monomial curves $(t^a, t^b, t^c)$
for pairwise coprime integers $a$, $b$, $c$ such that
$(a,b,c) \neq (1,1,1)$.
If such a ring is not finitely generated over a base field,
then it is a counterexample to the Hilbert's fourteenth problem.
Finite generation of such rings is deeply related to existence
of negative curves on certain normal projective surfaces.
We study a sufficient condition (Definition 3.6) for existence of a negative curve. Using it, we prove that, in the case of $(a+b+c)^2 > abc$, a negative curve exists.
Using a computer, we shall show that there exist examples in which this sufficient condition is not satisfied.}

\section{Introduction}
Let $k$ be a field.
Let $R$ be a polynomial ring over $k$ with finitely many variables.
For a field $L$ satisfying $k \subset L \subset Q(R)$,
Hilbert asked in 1900 whether the ring $L \cap R$ is finitely generated
as a $k$-algebra or not.
It is called the {\em Hilbert's fourteenth problem}.

The first counterexample to this problem 
was discovered by Nagata~\cite{N} in 1958.
An easier counterexample was found by Paul~C.~Roberts~\cite{Ro} in 1990.
Further counterexamples were given by Kuroda, Mukai, etc.

The Hilbert's fourteenth problem is deeply related to 
the following question of Cowsik~\cite{Cow}.
Let $R$ be a regular local ring (or a polynomial ring over a field).
Let $P$ be a prime ideal of $R$.
Cowsik asked whether the symbolic Rees ring 
\[
R_s(P) = \bigoplus_{r \geq 0} P^{(r)}T^r
\]
of $P$ is a
Noetherian ring or not.
His aim is to give a new approach to the Kronecker's problem,
that asks whether affine algebraic curves are set theoretic
complete intersection or not.
Kronecker's problem is still open, however, 
Roberts~\cite{R1} gave a counterexample to Cowsik's question in 1985.
Roberts constructed a regular local ring and a prime ideal
such that the completion coincides with Nagata's counterexample 
to the Hilbert's fourteenth problem.
In Roberts' example, the regular local ring contains a field
of characteristic zero, and the prime ideal splits after 
completion.
Later, Roberts~\cite{Ro} gave a new easier counterexample to both
Hilbert's fourteenth problem and Cowsik's question.
In his new example, the prime ideal does not split after 
completion, however, the ring still contains a field
of characteristic zero. 
It was proved that analogous rings of characteristic positive
are finitely generated (\cite{K1}, \cite{K2}).

On the other hand, let ${\frak p}_k(a,b,c)$ be the defining ideal
of the space monomial curves $(t^a, t^b, t^c)$ in $k^3$.
Then, ${\frak p}_k(a,b,c)$ is generated by 
at most three binomials in $k[x,y,z]$.
The symbolic Rees rings are deeply studied by many authors.
Huneke~\cite{Hu} and Cutkosky~\cite{Cu} developed
criterions for finite generation of such rings.
In 1994, Goto, Nishida and Watanabe~\cite{GNW}
proved that $R_s({\frak p}_k(7n-3, (5n-2)n, 8n-3))$ is 
not finitely generated over $k$ if the characteristic of $k$ is zero, 
$n \geq 4$ and $n \not\equiv 0 \ (3)$.
In their proof of infinite generation, they proved
the finite generation of 
$R_s({\frak p}_k(7n-3, (5n-2)n, 8n-3))$ in the case where $k$ is
of characteristic positive.
Goto and Watanabe conjectured that, for any $a$, $b$ and $c$, 
$R_s({\frak p}_k(a, b, c))$ is always
finitely generated over $k$ if the characteristic of $k$ is positive.

On the other hand, Cutkosky~\cite{Cu} gave a geometric meaning to 
the symbolic Rees ring $R_s({\frak p}_k(a, b, c))$.
Let $X$ be the blow-up of the weighted projective space
$\proj(k[x,y,z])$ at the smooth point $V_+({\frak p}_k(a, b, c))$.
Let $E$ be the exceptional curve of the blow-up.
Finite generation of $R_s({\frak p}_k(a, b, c))$ 
is equivalent to that of the total coordinate ring 
\[
TC(X) = \bigoplus_{D \in {\rm Cl}(X)} H^0(X, {\cal O}_X(D))
\]
of $X$.
If $-K_X$ is ample, one can prove that $TC(X)$ is finitely generated
using the cone theorem (cf.\ \cite{KMM}) as in \cite{kee-git}.
Cutkosky proved that $TC(X)$ is finitely generated
if $(-K_X)^2 > 0$, or equivalently $(a+b+c)^2 > abc$.
Finite generation of $TC(X)$ is deeply related
to existence of a negative curve $C$, i.e., a curve $C$ 
on $X$ satisfying $C^2 < 0$ and $C \neq E$.
In fact, in the case where $\sqrt{abc} \not\in {\Bbb Z}$, 
a negative exists if $TC(X)$ is finitely generated.
If a negative exists in the case 
where the characteristic of $k$ is positive, then
$TC(X)$ is finitely generated  by a result of M.~Artin~\cite{A}.

By a standard method (mod $p$ reduction), 
if there exists a negative curve in the case of characteristic zero,
then one can prove that a negative curve exists 
in the case of characteristic positive, therefore, 
$R_s({\frak p}_k(a, b, c))$ is finitely generated 
in the case of characteristic positive (cf.\ Lemma~\ref{lem3.2}).
In the examples of Goto-Nishida-Watanabe~\cite{GNW},
a negative curve exists, however, 
$R_s({\frak p}_k(a, b, c))$ is not finitely generated 
over $k$ in the case where $k$ is of characteristic zero 
(cf.\ Remark~\ref{rem3.3} below).

In Section~\ref{sec2}, we shall prove  that if 
$R_s({\frak p}_k(a, b, c))$ is not finitely generated,
then it is a counterexample to the Hilbert's fourteenth problem
(cf.\ Theorem~\ref{th2.1} and Remark~\ref{rem2.2}).

In Section~\ref{sec3}, we review some basic facts 
on finite generation of $R_s({\frak p}_k(a, b, c))$.
We define sufficient conditions for $X$ to have a negative curve
(cf.\ Definition~\ref{C123}).

In Section~\ref{sec4}, we shall prove that
there exists a negative curve in the case where
$(a+b+c)^2 > abc$ (cf.\ Theorem~\ref{th4.3}).
We should mention that if $(a+b+c)^2 > abc$, then
Cutkosky~\cite{Cu} proved that 
$R_s({\frak p}_k(a, b, c))$ is finitely generated.
Moreover if we assume $\sqrt{abc} \not\in {\Bbb Z}$,
existence of a negative curve follows from finite generation.
Existence of negative curves in these cases is an immediate 
conclusion of the cone theorem.
Our proof of existence of a negative curve is very simple, purely algebraic,
and do not need the cone theorem as Cutkosky's proof.

In Section~\ref{sec5}, we discuss the degree of a negative curve 
(cf.\ Theorem~\ref{bound}).
It is used in a computer programming in Section~\ref{c2}.

In Section~6.1, we prove that there exist examples in which a sufficient condition ((C2) in Definition 3.6) is not satisfied using a computer. In Section~6.2, we give a computer programming to check whether a negative curve exist or not.

\section{Symbolic Rees rings of monomial curves and
Hilbert's fourteenth problem}\label{sec2}

Throughout of this paper, we assume that rings are commutative
with unit.

For a prime ideal $P$ of a ring $A$,
$P^{(r)}$ denotes the $r$-th symbolic power of $P$, i.e.,
\[
P^{(r)} = P^rA_P \cap A .
\]
By definition, it is easily seen that $P^{(r)}P^{(r')} \subset
P^{(r+r')}$ for any $r, \ r' \geq 0$,
therefore,
\[
\bigoplus_{r \geq 0}P^{(r)}T^r
\]
is a subring of the polynomial ring $A[T]$.
This subring is called the {\em symbolic Rees ring} of $P$,
and denoted by $R_s(P)$.

Let $k$ be a field and $m$ be a positive integer.
Let $a_1$, \ldots, $a_m$ be positive integers.
Consider the $k$-algebra homomorphism
\[
\phi_k : k[x_1, \ldots, x_m] \longrightarrow
k[t]
\]
given by $\phi_k(x_i) = t^{a_i}$ for $i = 1, \ldots, m$,
where $x_1$, \ldots, $x_m$, $t$ are indeterminates over $k$.
Let ${\frak p}_k(a_1, \ldots, a_m)$ be the kernel of $\phi_k$.
We sometimes denote ${\frak p}_k(a_1, \ldots, a_m)$ simply by
${\frak p}$ or ${\frak p}_k$ if no confusion is possible.

\begin{Theorem}\label{th2.1}
Let $k$ be a field and $m$ be a positive integer.
Let $a_1$, \ldots, $a_m$ be positive integers.
Consider the prime ideal ${\frak p}_k(a_1, \ldots, a_m)$ 
of the polynomial ring $k[x_1, \ldots, x_m]$.

Let $\alpha_1$, $\alpha_2$, $\beta_1$, \ldots, $\beta_m$, $t$, $T$ 
be indeterminates over $k$.
Consider the following injective $k$-homomorphism
\[
\xi : k[x_1, \ldots, x_m, T] \longrightarrow
k(\alpha_1, \alpha_2, \beta_1, \ldots, \beta_m, t)
\]
given by $\xi(T) = \alpha_2/\alpha_1$ and $\xi(x_i) = \alpha_1\beta_i + t^{a_i}$ 
for $i = 1, \ldots, m$.

Then,
\[
k(\alpha_1\beta_1 + t^{a_1}, \alpha_1\beta_2 + t^{a_2}, \ldots, \alpha_1\beta_m + t^{a_m}, \alpha_2/\alpha_1) \cap k[\alpha_1, \alpha_2, \beta_1, \ldots, \beta_m, t]
= \xi(R_s({\frak p}_k(a_1, \ldots, a_m))) 
\]
holds true.
\end{Theorem}

\proof
Set $L = k(\alpha_1\beta_1 + t^{a_1}, \ldots, \alpha_1\beta_m + t^{a_m}, \alpha_2/\alpha_1)$.
Set $d = GCD(a_1, \ldots, a_m)$.
Then, $L$ is included in 
$k(\alpha_1, \alpha_2, \beta_1, \ldots, \beta_m, t^d)$.
Since
\[
k[\alpha_1, \alpha_2, \beta_1, \ldots, \beta_m, t] 
\cap k(\alpha_1, \alpha_2, \beta_1, \ldots, \beta_m, t^d)
= k[\alpha_1, \alpha_2, \beta_1, \ldots, \beta_m, t^d] ,
\]
we obtain the equality
\[
L \cap k[\alpha_1, \alpha_2, \beta_1, \ldots, \beta_m, t] =
L \cap k[\alpha_1, \alpha_2, \beta_1, \ldots, \beta_m, t^d]  .
\]
By the commutativity of the diagram
\[
\begin{array}{ccccc}
& & L & & \\
& & \downarrow & & \\
k[x_1, \ldots, x_m, T] & \longrightarrow & 
k(\alpha_1, \alpha_2, \beta_1, \ldots, \beta_m, t^d)
& \supset & k[\alpha_1, \alpha_2, \beta_1, \ldots, \beta_m, t^d]
\\
& {\scriptstyle \xi}\searrow \phantom{\scriptstyle \xi} &
\downarrow & & \downarrow \\
& & k(\alpha_1, \alpha_2, \beta_1, \ldots, \beta_m, t)
& \supset & k[\alpha_1, \alpha_2, \beta_1, \ldots, \beta_m, t]
\end{array}
\]
it is enough to prove this theorem in the case where
$GCD(a_1, \ldots, a_m) = 1$.

In the rest of this proof, we assume $GCD(a_1, \ldots, a_m) = 1$.

Consider the following injective $k$-homomorphism
\[
\tilde{\xi} : k[x_1, \ldots, x_m, T, t] \longrightarrow
k(\alpha_1, \alpha_2, \beta_1, \ldots, \beta_m, t)
\]
given by $\tilde{\xi}(T) = \alpha_2/\alpha_1$, $\tilde{\xi}(t) = t$ and 
$\tilde{\xi}(x_i) = \alpha_1\beta_i + t^{a_i}$ 
for $i = 1, \ldots, m$.
Here, remark that $\alpha_2/\alpha_1$, $\alpha_1\beta_1 + t^{a_1}$,
$\alpha_1\beta_2 + t^{a_2}$, \ldots, $\alpha_1\beta_m + t^{a_m}$, $t$
are algebraically independent over $k$.
By definition, the map $\xi$ is the restriction 
of $\tilde{\xi}$ to $k[x_1, \ldots, x_m, T]$.

We set $S = k[x_1, \ldots, x_m]$ and
$A = k[x_1, \ldots, x_m, t]$.
Let ${\frak q}$ be the ideal of $A$ generated by $x_1-t^{a_1}$, \ldots,
$x_m-t^{a_m}$.
Then ${\frak q}$ is the kernel of the map 
$\tilde{\phi_k} : A \rightarrow k[t]$
given by $\tilde{\phi_k}(t) = t$ and $\tilde{\phi_k}(x_i) = t^{a_i}$ for 
each $i$.
Since $\phi_k$ is the restriction of $\tilde{\phi_k}$ to $S$,
${\frak q} \cap S = {\frak p}$ holds.

Now we shall prove ${\frak q}^r \cap S = {\frak p}^{(r)}$ for each $r>0$.
Since ${\frak q}$ is a complete intersection,
${\frak q}^{(r)}$ coincides with ${\frak q}^r$ for any $r > 0$.
Therefore, 
it is easy to see ${\frak q}^r \cap S \supset {\frak p}^{(r)}$.

Since $GCD(a_1, \ldots, a_m) = 1$, there exists a monomial $M$ 
in $S$ 
such that $\phi_k(x_1^u) t = \phi_k(M)$ for some $u > 0$.
Let 
\[
\tilde{\phi_k}\otimes 1 : k[x_1, \ldots, x_m, x_1^{-1}, t] \longrightarrow 
k[t, t^{-1}]
\]
be the localization of $\tilde{\phi_k}$.
Then, the kernel of $\tilde{\phi_k}\otimes 1$ is equal to 
\[
{\frak q} k[x_1, \ldots, x_m, x_1^{-1}, t] =  
({\frak p}, t - \frac{M}{x_1^u}) k[x_1, \ldots, x_m, x_1^{-1}, t]
.
\]
Setting $t' = t - \frac{M}{x_1^u}$,
\[
{\frak q} A[x_1^{-1}] =
({\frak p}, t') k[x_1, \ldots, x_m, x_1^{-1}, t'] 
\]
holds.
Since $x_1$, \ldots, $x_m$, $t'$ are algebraically independent
over $k$,
\[
{\frak q}^r A[x_1^{-1}] \cap S[x_1^{-1}] =
({\frak p}, t')^r k[x_1, \ldots, x_m, x_1^{-1}, t'] \cap 
k[x_1, \ldots, x_m, x_1^{-1}]
= {\frak p}^r S[x_1^{-1}]
\]
holds.
Therefore,
\[
{\frak q}^r \cap S \subset
{\frak q}^r A[x_1^{-1}] \cap S = {\frak p}^r S[x_1^{-1}] \cap S
\subset  {\frak p}^{(r)} .
\]
We have completed the proof of ${\frak q}^r \cap S = {\frak p}^{(r)}$.

Let $R({\frak q})$ be the Rees ring of the ideal ${\frak q}$, i.e.,
\[
R({\frak q}) = \bigoplus_{r \geq 0}{\frak q}^rT^r \subset A[T] .
\]
Then, since ${\frak q}^r \cap S = {\frak p}^{(r)}$ for $r \geq 0$,
\[
R({\frak q}) \cap S[T] = R_s({\frak p})
\]
holds.
It is easy to verify
\[
R({\frak q}) \cap Q(S[T]) = R_s({\frak p}) 
\]
because $Q(S[T]) \cap A[T] = S[T]$,
where $Q( \ )$ means the field of fractions.
Here remark that $S[T] = k[x_1, \ldots, x_m, T]$ and
$A[T] = k[x_1, \ldots, x_m, T, t]$.
Therefore, we obtain the equality
\begin{equation}\label{eq1}
\tilde{\xi}(R({\frak q})) \cap L =
\xi(R_s({\frak p})) .
\end{equation}
Here, remember that $L$ is the field of fractions of ${\rm Im}(\xi)$.

On the other hand, setting $x'_i = x_i-t^{a_i}$ 
for $i = 1, \ldots, m$, we obtain the following:
\begin{eqnarray*}
R({\frak q}) & = & 
k[x_1, \ldots, x_m, x'_1T, \ldots, x'_mT, t] \\
& = & 
k[x'_1, \ldots, x'_m, x'_1T, \ldots, x'_mT, t] 
\end{eqnarray*}
Here, remark that $x'_1$, \ldots, $x'_m$, $T$, $t$ 
are algebraically independent over $k$.

By definition, $\tilde{\xi}(x'_i) = \alpha_1\beta_i$, 
and $\tilde{\xi}(x'_iT) = \alpha_2\beta_i$ for each $i$.

We set 
\begin{equation}\label{eq2}
B = \tilde{\xi}(R({\frak q}))
\end{equation}
and $C = k[\alpha_1, \alpha_2, \beta_1, \ldots, \beta_m, t]$.
Here, 
\[
B = 
\left(
k[\alpha_i\beta_j \mid i = 1, 2; \ j = 1, \ldots, m]
\right)
[t] \subset C .
\]
Since $B$ is a direct summand of $C$ as a $B$-module,
the equality
\begin{equation}\label{eq3}
C \cap Q(B) = B
\end{equation}
holds in $Q(C)$.

Then, since $L \subset Q(B)$,
we obtain 
\[
C \cap L
= (C \cap Q(B)) \cap L
= B \cap L
= \xi(R_s({\frak p})) 
\]
by the equations (\ref{eq1}), (\ref{eq2}) and (\ref{eq3}).
\qed

\begin{Remark}\label{rem2.2}
\begin{rm}
Let $k$ be a field.
Let $R$ be a polynomial ring over $k$ with finitely many variables.
For a field $L$ satisfying $k \subset L \subset Q(R)$,
Hilbert asked in 1900 whether the ring $L \cap R$ is finitely generated
as a $k$-algebra or not.
It is called the {\em Hilbert's fourteenth problem}.

The first counterexample to this problem 
was discovered by Nagata~\cite{N} in 1958.
An easier counterexample was found by Paul~C.~Roberts~\cite{Ro} in 1990.
Further counterexamples were given by Kuroda, Mukai, etc.

On the other hand, Goto, Nishida and Watanabe~\cite{GNW}
proved that $R_s({\frak p}_k(7n-3, (5n-2)n, 8n-3))$ is not finitely generated
over $k$ if the characteristic of $k$ is zero, $n \geq 4$ and 
$n \not\equiv 0 \ (3)$.
By Theorem~\ref{th2.1}, we know that they are new counterexamples
to the Hilbert's fourteenth problem.
\end{rm}
\end{Remark}

\begin{Remark}\label{rem2.3}
\begin{rm}
With notation as in Theorem~\ref{th2.1},
we set
\begin{eqnarray*}
D_1 & = &
\alpha_1 \frac{\partial \ \ }{\partial \alpha_1} 
+ \alpha_2 \frac{\partial \ \ }{\partial \alpha_2} 
- \beta_1 \frac{\partial \ \ }{\partial \beta_1} - \cdots - 
\beta_m \frac{\partial \ \ }{\partial \beta_m} \\
D_2 & = & 
a_1t^{a_1-1} \frac{\partial \ \ }{\partial \beta_1} + \cdots +
a_mt^{a_m-1} \frac{\partial \ \ }{\partial \beta_m} 
- \alpha_1 \frac{\partial \ \ }{\partial t} .
\end{eqnarray*}
Assume that the characteristic of $k$ is zero.

Then, one can prove that $\xi(R_s({\frak p}_k(a_1, \ldots, a_m)))$ 
is equal to the kernel of the derivations $D_1$ and $D_2$, i.e.,
\[
\xi(R_s({\frak p}_k(a_1, \ldots, a_m)))
= \{ f \in k[\alpha_1, \alpha_2, \beta_1, \ldots, \beta_m, t] \mid
D_1(f) = D_2(f) = 0 \} .
\]
\end{rm}
\end{Remark}

\section{Symbolic Rees rings of space monomial curves}\label{sec3}

In the rest of this paper, we restrict ourselves to the case $m = 3$.
For the simplicity of notation, we write $x$, $y$, $z$, $a$, $b$, $c$
for $x_1$, $x_2$, $x_3$, $a_1$, $a_2$, $a_3$, respectively.
We regard the polynomial ring $k[x, y, z]$ as a ${\Bbb Z}$-graded ring
by $\deg(x) = a$, $\deg(y) = b$ and $\deg(z) = c$.

${\frak p}_k(a, b, c)$ is the kernel of the $k$-algebra
homomorphism
\[
\phi_k : k[x,y,z] \longrightarrow k[t]
\]
given by $\phi_k(x) = t^a$, $\phi_k(y) = t^b$, $\phi_k(z) = t^c$.

By a result of Herzog~\cite{Her}, we know that ${\frak p}_k(a, b, c)$
is generated by at most three elements.
For example, ${\frak p}_k(3, 4, 5)$ is minimally generated by
$x^3-yz$, $y^2-zx$ and $z^2-x^2y$.
On the other hand, ${\frak p}_k(3, 5, 8)$ is minimally generated by
$x^5-y^3$ and $z-xy$.
We can choose a generating system of ${\frak p}_k(a, b, c)$ which is 
independent of $k$.

We are interested in the symbolic powers of ${\frak p}_k(a, b, c)$.
If ${\frak p}_k(a, b, c)$ is generated by two elements,
then the symbolic powers always coincide with ordinary powers
because ${\frak p}_k(a,b,c)$ is a complete intersection.
However, it is known that, if ${\frak p}_k(a, b, c)$ is minimally
generated by three elements, the second symbolic power is strictly
bigger than the second ordinary power.
For example, the element
\[
h = (x^3-yz)^2 - (y^2-zx)(z^2-x^2y)
\]
is contained in ${\frak p}_k(3, 4, 5)^2$, and
is divisible by $x$.
Therefore, $h/x$ is an element in ${\frak p}_k(3, 4, 5)^{(2)}$
of degree $15$.
Since $[{\frak p}_k(3, 4, 5)^2]_{15} = 0$, 
$h/x$ is not contained in ${\frak p}_k(3, 4, 5)^2$.

We are interested in finite generation of the symbolic Rees ring
$R_s({\frak p}_k(a, b, c))$.
It is known that this problem 
is reduced to the case where $a$, $b$ and $c$ are pairwise coprime, i.e.,
\[
(a,b) = (b,c) = (c,a) = 1 .
\]

In the rest of this paper, we always assume that 
$a$, $b$ and $c$ are pairwise coprime.

Let ${\Bbb P}_k(a,b,c)$ be the weighted projective space
$\proj(k[x,y,z])$. 
Then
\[
{\Bbb P}_k(a,b,c) \setminus
\{ V_+(x,y), \ \ V_+(y,z), \ \ V_+(z,x) \}
\]
is a regular scheme.
In particular, ${\Bbb P}_k(a,b,c)$ is smooth at the point 
$V_+({\frak p}_k(a, b, c))$.
Let $\pi : X_k(a,b,c) \rightarrow {\Bbb P}_k(a,b,c)$ be the blow-up 
at $V_+({\frak p}_k(a, b, c))$.
Let $E$ be the exceptional divisor, i.e.,
\[
E = \pi^{-1}(V_+({\frak p}_k(a, b, c))) .
\]

We sometimes denote ${\frak p}_k(a, b, c)$ 
(resp.\ ${\Bbb P}_k(a,b,c)$, $X_k(a,b,c)$ )
simply by ${\frak p}$ or ${\frak p}_k$
(resp.\ ${\Bbb P}$ or ${\Bbb P}_k$, 
$X$ or $X_k$)
if no confusion is possible.

It is easy to see that
\[
{\rm Cl}({\Bbb P}) = {\Bbb Z} H \simeq {\Bbb Z},
\]
where $H$ is a Weil divisor corresponding to the reflexive sheaf
${\cal O}_{\Bbb P}(1)$.
Set $H = \sum_im_iD_i$, where $D_i$'s are subvarieties of ${\Bbb P}$
of codimension one.
We may choose $D_i$'s such that $D_i \not\ni V_+({\frak p})$ for any $i$.
Then, set $A = \sum_im_i\pi^{-1}(D_i)$.

One can prove that
\[
{\rm Cl}(X) = {\Bbb Z}A + {\Bbb Z}E \simeq {\Bbb Z}^2 .
\]
Since all Weil divisors on $X$ are ${\Bbb Q}$-Cartier,
we have the intersection pairing
\[
{\rm Cl}(X) \times {\rm Cl}(X) \longrightarrow {\Bbb Q} ,
\]
that satisfies
\[
A^2 = \frac{1}{abc}, \ \ E^2 = -1, \ \ A.E = 0 .
\]
Therefore, we have
\[
(n_1A - r_1E).(n_2A - r_2E) = \frac{n_1n_2}{abc} - r_1r_2 .
\]
Here, we have the following natural identification:
\[
H^0(X, {\cal O}_X(nA - rE))
=
\left\{
\begin{array}{cl}
\left[
{\frak p}^{(r)}
\right]_n \ \ 
& (r \geq 0) \\
S_n & (r < 0)
\end{array}
\right.
\]
Therefore, the {\em total coordinate ring} (or {\em Cox ring})
\[
TC(X) = \bigoplus_{n, r \in {\Bbb Z}} H^0(X, {\cal O}_X(nA - rE))
\]
is isomorphic to the extended symbolic Rees ring
\[
R_s({\frak p})[T^{-1}] = \cdots \oplus ST^{-2} \oplus ST^{-1} \oplus
S \oplus {\frak p}T \oplus {\frak p}^{(2)}T^2 \oplus \cdots  .
\]
We refer the reader to Hu-Keel~\cite{kee-git}
for finite generation of total coordinate rings.
It is well-known that $R_s({\frak p})[T^{-1}]$ 
is Noetherian if and only if so is $R_s({\frak p})$.

\begin{Remark}\label{cri}
\begin{rm}
By Huneke's criterion~\cite{Hu} and a result of Cutkosky~\cite{Cu}, 
the following four
conditions are equivalent:
\begin{itemize}
\item[(1)]
$R_s({\frak p})$ is a Noetherian ring, or equivalently, 
finitely generated over $k$.
\item[(2)]
$TC(X)$ is a Noetherian ring, or equivalently, 
finitely generated over $k$.
\item[(3)]
There exist positive integers $r$, $s$, 
$f \in {\frak p}^{(r)}$, $g \in {\frak p}^{(s)}$, 
and $h \in (x,y,z) \setminus {\frak p}$ such that
\[
\ell_{S_{(x,y,z)}}(S_{(x,y,z)}/(f,g,h))
= rs \cdot \ell_{S_{(x,y,z)}}(S_{(x,y,z)}/({\frak p}, h)) ,
\]
where $\ell_{S_{(x,y,z)}}$ is the length as an 
$S_{(x,y,z)}$-module.
\item[(4)]
There exist curves $C$ and $D$ on $X$ such that
\[
C \neq D, \ \ C \neq E, \ \ D \neq E, \ \ C.D = 0 .
\]
\end{itemize}
Here, a curve means a closed irreducible reduced subvariety
of dimension one.

The condition (4) as above is equivalent to
that just one of the following two conditions is satisfied:
\begin{itemize}
\item[(4-1)]
There exist curves $C$ and $D$ on $X$ such that
\[
C \neq E, \ \ D \neq E, \ \ C^2 < 0, \ \ D^2 > 0, \ \ C.D = 0 .
\]
\item[(4-2)]
There exist curves $C$ and $D$ on $X$ such that
\[
C \neq E, \ \ D \neq E, \ \ C \neq D, \ \ C^2=D^2=0 .
\]
\end{itemize}
\end{rm}
\end{Remark}

\begin{Definition}
\begin{rm}
A curve $C$ on $X$ is called a {\em negative curve}
if 
\[
C \neq E \ \ \mbox{and} \ \ C^2 < 0  .
\]
\end{rm}
\end{Definition}

\begin{Remark}
\begin{rm}
Suppose that a divisor $F$ is linearly equivalent to $nA - rE$.
Then, we have 
\[
F^2 = \frac{{n}^2}{abc} - {r}^2 .
\]

If (4-2) in Remark~\ref{cri} is satisfied,
then all of $a$, $b$ and $c$ must be squares of integers
because $a$, $b$, $c$ are pairwise coprime.
In the case where one of $a$, $b$ and $c$ is not square,
the condition (4) is equivalent to (4-1).
Therefore, in this case, a negative curve exists
if $R_s({\frak p})$ is finitely generated over $k$.

Suppose $(a,b,c) =(1,1,1)$.
Then ${\frak p}$ coincides with $(x-y, y-z)$.
Of course, $R_s({\frak p})$ is a Noetherian ring since
the symbolic powers coincide with the ordinary powers.
By definition it is easy to see that there is no negative curve 
in this case,
therefore, (4-2) in Remark~\ref{cri} happens.

The authors know no other examples in which (4-2) happens.
\end{rm}
\end{Remark}

In the case of $(a,b,c) = (3,4,5)$,
the proper transform of
\[
V_+( \frac{(x^3-yz)^2 - (y^2-zx)(z^2-x^2y)}{x} )
\]
is the negative curve on $X$,
that is linearly equivalent to $15A - 2E$.

It is proved that two distinct negative curves never exist.

In the case where the characteristic of $k$ is positive,
Cutkosky~\cite{Cu} proved that  
$R_s({\frak p})$ is finitely generated over $k$
if there exists a negative curve on $X$.

We remark that there exists a negative curve on $X$
if and only if there exists positive integers $n$ and $r$ such that
\[
\frac{n}{r} < \sqrt{abc} \ \ \mbox{and} \ \ [{\frak p}^{(r)}]_n \neq 0  .
\]

We are interested in existence of a negative curve.
Let $a$, $b$ and $c$ be pairwise coprime positive integers. 
By the following lemma,
if there exists a negative curve on $X_{k_0}(a,b,c)$ for a field 
$k_0$ of characteristic $0$, 
then there exists a negative curve on $X_k(a,b,c)$ for any field $k$.

\begin{Lemma}\label{lem3.2}
Let $a$, $b$ and $c$ be pairwise coprime positive integers. 
\begin{enumerate}
\item
Let $K/k$ be a field extension.
Then, for any integers $n$ and $r$,
\[
[{\frak p}_k(a,b,c)^{(r)}]_n \otimes_kK =
[{\frak p}_K(a,b,c)^{(r)}]_n .
\]
\item
For any integers $n$, $r$ and any prime number $p$,
\[
\dim_{{\Bbb F}_p} [{\frak p}_{{\Bbb F}_p}(a,b,c)^{(r)}]_n \geq
\dim_{{\Bbb Q}} [{\frak p}_{{\Bbb Q}}(a,b,c)^{(r)}]_n 
\]
holds, where ${\Bbb Q}$ is the field of rational numbers,
and ${\Bbb F}_p$ is the prime field of characteristic $p$.
Here, $\dim_{{\Bbb F}_p}$ (resp.\ $\dim_{{\Bbb Q}}$) denotes
the dimension as an ${\Bbb F}_p$-vector space
(resp.\ ${\Bbb Q}$-vector space).
\end{enumerate}
\end{Lemma}

\proof
Since $S \rightarrow S \otimes_kK$ is flat, it is easy to prove 
the assertion (1).

We shall prove the assertion (2).
Let ${\Bbb Z}$ be the ring of rational integers.
Set $S_{\Bbb Z} = {\Bbb Z}[x,y,z]$.
Let ${\frak p}_{\Bbb Z}$ be 
the kernel of the ring homomorphism
\[
\phi_{\Bbb Z} : S_{\Bbb Z} \longrightarrow
{\Bbb Z}[t]
\]
given by $\phi_{\Bbb Z}(x) = t^{a}$,
$\phi_{\Bbb Z}(y) = t^{b}$ and
$\phi_{\Bbb Z}(z) = t^{c}$.
Since the cokernel of $\phi_{\Bbb Z}$ is 
${\Bbb Z}$-free module, we know
\[
{\frak p}_{\Bbb Z} \otimes_{\Bbb Z}k
= \ker(\phi_{\Bbb Z}) \otimes_{\Bbb Z}k
= \ker(\phi_k) = 
{\frak p}_k
\]
for any field $k$.

Consider the following exact sequence of
${\Bbb Z}$-free modules:
\[
0 \longrightarrow {{\frak p}_{\Bbb Z}}^{(r)}
\longrightarrow S_{\Bbb Z}
\longrightarrow S_{\Bbb Z}/{{\frak p}_{\Bbb Z}}^{(r)}
\longrightarrow 0 
\]
For any field $k$,
the following sequence is exact:
\[
0 \longrightarrow {{\frak p}_{\Bbb Z}}^{(r)} \otimes_{\Bbb Z}k
\longrightarrow S
\longrightarrow S_{\Bbb Z}/{{\frak p}_{\Bbb Z}}^{(r)} \otimes_{\Bbb Z}k
\longrightarrow 0 
\]

Since ${\frak p}_{\Bbb Z}S_{\Bbb Z}[x^{-1}]$ is generated by a regular
sequence, we know
\[
{{\frak p}_{\Bbb Z}}^{(r)}S_{\Bbb Z}[x^{-1}]
= {{\frak p}_{\Bbb Z}}^rS_{\Bbb Z}[x^{-1}]
\]
for any $r \geq 0$.
Therefore, for any $f \in {{\frak p}_{\Bbb Z}}^{(r)}$,
there is a positive integer $u$ such that
\[
x^u f \in {{\frak p}_{\Bbb Z}}^r .
\]

Let $p$ be a prime number.
Consider the natural surjective ring homomorphism
\[
\eta : S_{\Bbb Z}
\longrightarrow
S_{\Bbb Z} \otimes_{\Bbb Z} {\Bbb F}_p .
\]
Suppose $f \in {{\frak p}_{\Bbb Z}}^{(r)}$.
Since $x^u f \in {{\frak p}_{\Bbb Z}}^r$ for some positive integer
$u$, we obtain
\[
x^u \eta(f) \in \eta({{\frak p}_{\Bbb Z}}^r)
= {{\frak p}_{{\Bbb F}_p}}^r .
\]
Hence we know
\[
{{\frak p}_{\Bbb Z}}^{(r)}\otimes_{\Bbb Z} {\Bbb F}_p 
= \eta({{\frak p}_{\Bbb Z}}^{(r)})
\subset 
{{\frak p}_{{\Bbb F}_p}}^{(r)} .
\]
We obtain
\[
\rank_{{\Bbb Z}} [{{\frak p}_{\Bbb Z}}^{(r)}]_n
= \dim_{{\Bbb F}_p} [{{\frak p}_{\Bbb Z}}^{(r)}]_n \otimes_{\Bbb Z} {\Bbb F}_p
\leq \dim_{{\Bbb F}_p} [{{\frak p}_{{\Bbb F}_p}}^{(r)}]_n 
\]
for any $r \geq 0$ and $n \geq 0$.
Here, $\rank_{{\Bbb Z}}$ denotes the rank as a ${\Bbb Z}$-module.

On the other hand, it is easy to see that
\[
{{\frak p}_{\Bbb Z}}^{(r)} \otimes_{\Bbb Z}{\Bbb Q}
= {{\frak p}_{\Bbb Q}}^{(r)}
\]
for any $r \geq 0$.
Therefore, we have
\[
\rank_{\Bbb Z} [{{\frak p}_{\Bbb Z}}^{(r)}]_n
= \dim_{\Bbb Q} [{{\frak p}_{\Bbb Q}}^{(r)}]_n
\]
for any $r \geq 0$ and $n \geq 0$.

Hence, we obtain
\[
\dim_{\Bbb Q} [{{\frak p}_{\Bbb Q}}^{(r)}]_n
\leq 
\dim_{{\Bbb F}_p} [{{\frak p}_{{\Bbb F}_p}}^{(r)}]_n
\]
for any $r \geq 0$, $n \geq 0$, and any prime number $p$.
\qed

\begin{Remark}\label{rem3.3}
\begin{rm}
Let $a$, $b$, $c$ be pairwise coprime positive integers.
Assume that there exists a negative curve on $X_{k_0}(a,b,c)$
for a field ${k_0}$ of characteristic zero.

By Lemma~\ref{lem3.2}, we know that
there exists a negative curve on $X_k(a,b,c)$
for any field $k$.
Therefore, if $k$ is a field of characteristic positive,
then the symbolic Rees ring $R_s({\frak p}_k)$ is finitely
generated over $k$ by a result of Cutkosky~\cite{Cu}.
However, if $k$ is a field of characteristic zero,
then $R_s({\frak p}_k)$ is not necessary Noetherian.
In fact, assume that $k$ is of characteristic zero
and $(a,b,c) = (7n-3,(5n-2)n,8n-3)$ with $n \not\equiv 0 \ (3)$ and $n \geq 4$
as in Goto-Nishida-Watanabe~\cite{GNW}.
Then there exists a negative curve, but 
$R_s({\frak p}_k)$ is not Noetherian.
\end{rm}
\end{Remark}

\begin{Definition}\label{C123}
\begin{rm}
Let $a$, $b$, $c$ be pairwise coprime positive integers.
Let $k$ be a field.

We define the following three conditions:
\begin{itemize}
\item[(C1)]
There exists a negative curve on $X_k(a,b,c)$, i.e.,
$[{{\frak p}_k(a,b,c)}^{(r)}]_n \neq 0$ for some positive integers $n$, $r$ 
satisfying $n/r < \sqrt{abc}$.
\item[(C2)]
There exist positive integers $n$, $r$ 
satisfying $n/r < \sqrt{abc}$ and $\dim_k S_n > r(r+1)/2$.
\item[(C3)]
There exist positive integers $q$, $r$ 
satisfying $abcq/r < \sqrt{abc}$ and $\dim_k S_{abcq} > r(r+1)/2$.
\end{itemize}
Here, $\dim_k$ denotes the dimension as a $k$-vector space.
\end{rm}
\end{Definition}

By the following lemma, we know the implications
\[
(C3) \Longrightarrow (C2) \Longrightarrow (C1)
\]
since $\dim_k [{\frak p}^{(r)}]_n = \dim_k S_n - \dim_k [S/{\frak p}^{(r)}]_n$.

\begin{Lemma}
Let $a$, $b$, $c$ be pairwise coprime positive integers.
Let $r$ and $n$ be non-negative integers.
Then, 
\[
\dim_k [S/{\frak p}^{(r)}]_n \leq r(r+1)/2 
\]
holds true for any field $k$.
\end{Lemma}

\proof
Since $x$, $y$, $z$ are non-zero divisors on $S/{\frak p}^{(r)}$,
we have only to prove that
\[
\dim_k [S/{\frak p}^{(r)}]_{abcq} = r(r+1)/2
\]
for $q \gg 0$.

The left-hand side is the multiplicity of the $abc$-th
Veronese subring
\[
[S/{\frak p}^{(r)}]^{(abc)}
= \oplus_{q \geq 0}
[S/{\frak p}^{(r)}]_{abcq} .
\]
Therefore, for $q \gg 0$, we have
\begin{eqnarray*}
\dim_k [S/{\frak p}^{(r)}]_{abcq} & = & 
\ell ([S/{\frak p}^{(r)} + (x^{bc})]^{(abc)}) \\
& = & e((x^{bc}), [S/{\frak p}^{(r)}]^{(abc)}) \\
& = & \frac{1}{abc} e((x^{bc}), S/{\frak p}^{(r)}) \\
& = & \frac{1}{a} e((x), S/{\frak p}^{(r)}) \\
& = & \frac{1}{a} e((x), S/{\frak p})
\ell_{S_{\frak p}}(S_{\frak p}/{\frak p}^rS_{\frak p}) \\
& = & \frac{r(r+1)}{2}
\end{eqnarray*}
\qed

\begin{Remark}\label{r=01}
\begin{rm}
It is easy to see that $[{\frak p}_k(a,b,c)]_n \neq 0$
if and only if $\dim_k S_n \geq 2$.
Therefore, if we restrict ourselves to $r = 1$, 
then (C1) and (C2) are equivalent.

However, even if $[{\frak p}_k(a,b,c)^{(2)}]_n \neq 0$,
$\dim_k S_n$ is not necessary bigger than $3$.
In fact, since ${\frak p}_k(5,6,7)$ contains $y^2-zx$,
we know $[{\frak p}_k(5,6,7)^2]_{24} \neq 0$.
In this case, $\dim_k S_{24}$ is equal to three.

Here assume that (C1) is satisfied for $r = 2$.
Furthermore, we assume that the characteristic of $k$ is zero.
Then, there exists $f \neq 0$ in 
$[{\frak p}_k(a,b,c)^{(2)}]_n$ such that $n < 2 \sqrt{abc}$
for some $n > 0$.
Let $f = f_1\cdots f_s$ be the irreducible decomposition.
Then, at least one of $f_i$'s satisfies the condition (C1).
If it satisfies (C1) with $r = 1$, then (C2) is satisfied as above.
Suppose that the irreducible component satisfies (C2) with $r = 2$.
For the simplicity of notation,
we assume that $f$ itself is irreducible.
We want to show $\dim_k S_n \geq 4$.
Assume the contrary.
By Lemma~\ref{lem3.2}~(1), we may assume that $f$ is a polynomial
with rational coefficients.
Set
\[
f = k_1 x^{\alpha_1}y^{\beta_1}z^{\gamma_1} -
k_2 x^{\alpha_2}y^{\beta_2}z^{\gamma_2} + 
k_3 x^{\alpha_3}y^{\beta_3}z^{\gamma_3} .
\]
Furthermore, we may assume that $k_1$, $k_2$, $k_3$ are non-negative integers
such that $GCD(k_1,k_2,k_3) = 1$.
Since
\[
\frac{\partial f}{\partial x}, 
\frac{\partial f}{\partial y}, 
\frac{\partial f}{\partial z} \in {\frak p}_k(a,b,c) 
\]
as in Remark~\ref{4.1},
we have
\[
\left(
\begin{array}{ccc}
\alpha_1 & \alpha_2 & \alpha_3 \\
\beta_1 & \beta_2 & \beta_3 \\
\gamma_1 & \gamma_2 & \gamma_3
\end{array}
\right)
\left(
\begin{array}{r}
k_1 \\ -k_2 \\ k_3 
\end{array}
\right)
 =
\left(
\begin{array}{c}
0 \\ 0 \\ 0
\end{array}
\right) .
\]
Therefore, we have
\[
(x^{\alpha_1}y^{\beta_1}z^{\gamma_1})^{k_1}
(x^{\alpha_3}y^{\beta_3}z^{\gamma_3})^{k_3}
= 
(x^{\alpha_2}y^{\beta_2}z^{\gamma_2})^{k_2} .
\]
Since $f$ is irreducible, $x^{\alpha_1}y^{\beta_1}z^{\gamma_1}$ and 
$x^{\alpha_3}y^{\beta_3}z^{\gamma_3}$ have no common divisor.
Note that $k_2 = k_1+k_3$ since $f \in {\frak p}_k(a,b,c)$.
Since $k_1$ and $k_3$ are relatively prime,
there exist monomials $N_1$ and $N_3$ such that
$x^{\alpha_1}y^{\beta_1}z^{\gamma_1} = N_1^{k_1+k_3}$, 
$x^{\alpha_3}y^{\beta_3}z^{\gamma_3} = N_3^{k_1+k_3}$
and
$x^{\alpha_2}y^{\beta_2}z^{\gamma_2} = N_1^{k_1}N_3^{k_3}$.
Then
\[
f = k_1 N_1^{k_1+k_3} - (k_1+k_3)N_1^{k_1}N_3^{k_3} + k_3 N_3^{k_1+k_3} .
\]
Then, $f$ is divisible by $N_1-N_3$.
Since $f$ is irreducible, $f$ is equal to $N_1-N_3$.
It contradicts to 
\[
\frac{\partial f}{\partial x}, 
\frac{\partial f}{\partial y}, 
\frac{\partial f}{\partial z} \in {\frak p}_k(a,b,c) .
\]

Consequently, if (C1) is satisfied with $r \leq 2$ for a field $k$
of characteristic zero,
then (C2) is satisfied.

We shall discuss the difference between (C1) and (C2) in 
Section~\ref{c2}.
\end{rm}
\end{Remark}

\begin{Remark}\label{ci}
\begin{rm}
Let $a$, $b$ and $c$ be pairwise coprime positive integers.
Assume that ${\frak p}_k(a,b,c)$ is a complete intersection,
i.e., generated by two elements.

Permuting $a$, $b$ and $c$, we may assume that
\[
{\frak p}_k(a,b,c)
= (
x^b-y^a, \ z - x^\alpha y^\beta
) 
\]
for some $\alpha, \beta \geq 0$ satisfying $\alpha a + \beta b = c$.
If $ab < c$, then 
\[
\deg(x^b-y^a) = ab < \sqrt{abc} .
\]
If $ab > c$, then 
\[
\deg(z - x^\alpha y^\beta) = c < \sqrt{abc} .
\]
If $ab = c$, then $(a,b,c)$ must be equal to $(1,1,1)$.
Ultimately, there exists a negative curve if $(a,b,c) \neq (1,1,1)$.
\end{rm}
\end{Remark}

\section{The case where $(a+b+c)^2 > abc$}\label{sec4}

In the rest of this paper, we set $\xi = abc$ and $\eta = a+b+c$
for pairwise coprime positive integers $a$, $b$ and $c$.

For $v = 0, 1, \ldots, \xi-1$,
we set
\[
S^{(\xi,v)} = \oplus_{q \geq 0} S_{\xi q + v} .
\]
This is a module over $S^{(\xi)} = \oplus_{q \geq 0} S_{\xi q}$.

\begin{Lemma}\label{dim}
\[
\dim_k [S^{(\xi,v)}]_q = \dim_k S_{\xi q + v}
= \frac{1}{2} \left\{
\xi q^2 + (\eta + 2v) q + 2 \dim_k S_v \right\}
\]
holds for any $q \geq 0$.
\end{Lemma}

The following simple proof is due to Professor Kei-ichi Watanabe.
We appreciate him very much.

\proof
We set $a_n = \dim_k S_n$ for each integer $n$.
Set
\[
f(t) = \sum_{n \in {\Bbb Z}}a_nt^n .
\]
Here we put $a_n = 0$ for $n < 0$.
Then, the equality
\[
f(t) = \frac{1}{(1-t^a)(1-t^b)(1-t^c)}
\]
holds.

Set $b_n = a_n - a_{n - \xi}$.
Then, $b_n$ is equal
to the coefficient of $t^n$ in $(1-t^\xi)f(t)$ for each $n$.
Furthermore, $b_n - b_{n-1}$ 
is equal to the coefficient of $t^n$ in $(1-t)(1-t^\xi)f(t)$ for each $n$.

On the other hand, we have the equality
\begin{equation}\label{*}
(1-t)(1-t^\xi)f(t)
= g(t) \times \frac{1}{1-t} 
= g(t) \times (1 + t + t^2 + \cdots) ,
\end{equation}
where
\[
g(t) = \frac{1 + t + \cdots + t^{\xi - 1}}{
(1 + t + \cdots + t^{a-1})(1 + t + \cdots + t^{b-1})
(1 + t + \cdots + t^{c-1})
} .
\]
Since $a$, $b$ and $c$ are pairwise coprime,
$g(t)$ is a polynomial of degree $\xi - \eta + 2$.
Therefore, the coefficient of $t^n$ in $(1-t)(1-t^\xi)f(t)$
is equal to $g(1)$ for $n \geq \xi - \eta + 2$
by the equation (\ref{*}).
It is easy to see $g(1) = 1$.

Since $b_n - b_{n-1} = 1$ for $n \geq \xi + 1$,
\[
b_n = b_\xi + (n - \xi)
\]
holds for any $n \geq \xi$.
Then, 
\begin{eqnarray*}
a_{\xi q + v} - a_v & = & 
\sum_{i = 1}^q (a_{\xi i + v} - a_{\xi(i-1) + v}) \\
& = & \sum_{i = 1}^q b_{\xi i + v} \\
& = & \sum_{i = 1}^q \left( b_\xi + \xi(i-1) + v
\right) \\
& = & b_\xi q + \xi \frac{(q-1)q}{2} + vq \\
& = & \frac{\xi}{2} q^2 + \left( 
b_\xi - \frac{\xi}{2} + v \right) q .
\end{eqnarray*}

Recall that $b_\xi$ is the coefficient of $t^\xi$ in 
\begin{equation}\label{eq10}
(1-t^\xi)f(t) = \frac{g(t)}{(1-t)^2}
= g(t) \times \left( 1 + 2t + \cdots + (n+1)t^n + \cdots \right) .
\end{equation}
Setting
\[
g(t) = c_0 + c_1t + \cdots + c_{\xi - \eta + 2}t^{\xi - \eta + 2} ,
\]
it is easy to see
\begin{equation}\label{eq11}
c_i = c_{\xi - \eta + 2 - i}
\end{equation}
for each $i$.
Therefore, by the equations (\ref{eq10}) and (\ref{eq11}), we have
\[
b_\xi = c_0 (\xi + 1) + c_1 \xi + \cdots + 
c_{\xi - \eta + 2} (\eta - 1)
= (c_0 + c_1 + \cdots + c_{\xi - \eta + 2}) \times 
\frac{\xi + \eta}{2} .
\]
Since $g(1) = 1$, we have $b_\xi = \frac{\xi + \eta}{2}$.
Thus, 
\[
a_{\xi q + v} = 
\frac{\xi}{2} q^2 + \left( 
\frac{\xi + \eta}{2} - \frac{\xi}{2} + v \right) q + a_v .
\]
\qed

Before proving Theorem~\ref{th4.3}, we need the following lemma:

\begin{Lemma}\label{012}
\begin{rm}
Assume that $a$, $b$ and $c$ are pairwise coprime
positive integers such that $(a,b,c) \neq (1,1,1)$.
Then, $\eta - \sqrt{\xi} \neq 0, 1, 2$.
\end{rm}
\end{Lemma}

\proof
We may assume that all of $a$, $b$ and $c$ are squares of integers.
It is sufficient to show
\[
\alpha^2 + \beta^2 + \gamma^2 - \alpha\beta\gamma \neq 0, 1, 2
\]
for pairwise coprime positive integers $\alpha$, $\beta$, $\gamma$
such that $(\alpha, \beta, \gamma) \neq (1,1,1)$.

Assume the contrary.
Suppose that $(\alpha_0, \beta_0, \gamma_0)$ is a counterexample
such that $\alpha_0 + \beta_0 + \gamma_0$ is minimum.
We may assume $1 \leq \alpha_0 \leq \beta_0 \leq \gamma_0$.

Set
\[
f(x) = x^2 - \alpha_0\beta_0x + \alpha_0^2 + \beta_0^2 .
\]

First suppose $\alpha_0\beta_0 \leq \gamma_0$.
Then,
\[
f(\gamma_0) \geq f(\alpha_0\beta_0) = \alpha_0^2 + \beta_0^2 
\geq 2 .
\]
Since $f(\gamma_0) = 0, 1$, or $2$, we have 
\[
\gamma_0 = \alpha_0 \beta_0 \ \
\mbox{and} \ \ 
\alpha_0^2 + \beta_0^2 = 2 .
\]
Then, we obtain the equality
 $\alpha_0 = \beta_0 = \gamma_0 = 1$ immediately.
It is a contradiction.

Next, suppose $\frac{\alpha_0\beta_0}{2} < \gamma_0 < \alpha_0\beta_0$.
Then, $0 < \alpha_0\beta_0 -\gamma_0 < \gamma_0$ and
\[
f(\alpha_0\beta_0 -\gamma_0) = f(\gamma_0) = 0, 1,\ \mbox{or} \  2 .
\]
It is easy to see that $\alpha_0$, $\beta_0$, $\alpha_0\beta_0 -\gamma_0$
are pairwise coprime positive integers.
By the minimality of $\alpha_0 + \beta_0 + \gamma_0$, 
we have $\alpha_0 = \beta_0 = \alpha_0\beta_0 -\gamma_0 = 1$.
Then, $\gamma_0$ must be zero.
It is a contradiction.

Finally, suppose $0 < \gamma_0 \leq \frac{\alpha_0\beta_0}{2}$.
Since $\beta_0 \leq \gamma_0 \leq \frac{\alpha_0\beta_0}{2}$,
we have $\alpha_0 \geq 2$.
If $\alpha_0 = 2$, then $2 \leq \beta_0 = \gamma_0$.
It contradicts to $(\beta_0, \gamma_0) = 1$.
Assume $\alpha_0 \geq 3$.
Since $\beta_0 < \gamma_0$,
\[
f(\gamma_0) < f(\beta_0) = (2 - \alpha_0) \beta_0^2 + \alpha_0^2
\leq 0 .
\]
It is a contradiction.
\qed

\begin{Theorem}\label{th4.3}
Let $a$, $b$ and $c$ be pairwise coprime integers such that
$(a, b, c) \neq (1,1,1)$.

Then, we have the following:
\begin{enumerate}
\item
Assume that $\sqrt{abc} \not\in {\Bbb Z}$.
Then, (C3) holds if and only if $(a+b+c)^2 > abc$.
\item
Assume that $\sqrt{abc} \in {\Bbb Z}$.
Then, (C3) holds if and only if $(a+b+c)^2 > 9abc$.
\item
If $(a+b+c)^2 > abc$, then, (C2) holds.
In particular, a negative curve exists in this case.
\end{enumerate}
\end{Theorem}

\proof
%
Remember that,
by Lemma~\ref{dim}, we obtain
\[
\dim_k S_{\xi q}
= \frac{1}{2} (
\xi q^2 + \eta q + 2  ) 
\]
for any $q \geq 0$.

First we shall prove the assertion (1).
Assume that (C3) is satisfied.
Then, 
\[
\left\{
\begin{array}{l}
\sqrt{\xi} > \frac{\xi q}{r} \\
\frac{\xi q^2 + \eta q + 2}{2} > \frac{r(r+1)}{2}
\end{array}
\right.
\]
is satisfied for some positive integers $r$ and $q$.
The second inequality is equivalent to
$\xi q^2 + \eta q \geq r(r+1)$ since both
integers are even.
Since
\[
\xi q^2 + \eta q \geq r^2 + r > \xi q^2 + \sqrt{\xi} q ,
\]
we have $\eta > \sqrt{\xi}$ immediately.

Assume $\eta > \sqrt{\xi}$ and $\sqrt{\xi} \not\in {\Bbb Z}$.
Let $\epsilon$ be a real number satisfying $0 < \epsilon < 1$ and
\begin{equation}\label{**}
2 \epsilon \sqrt{\xi} < \frac{\eta - \sqrt{\xi}}{2} .
\end{equation}
Since $\sqrt{\xi} \not\in {\Bbb Q}$, there exist positive integers
$r$ and $q$ such that
\[
\epsilon > r - \sqrt{\xi} q > 0 .
\]
Then, 
\[
\frac{r}{q} < \sqrt{\xi} + \frac{\epsilon}{q}
\leq \sqrt{\xi} + \epsilon <
\sqrt{\xi} + \frac{\eta - \sqrt{\xi}}{2}
=  \frac{\eta + \sqrt{\xi}}{2} .
\]
Since $\sqrt{\xi} q + \epsilon > r$, we have
\[
\xi q^2 + 2  \epsilon \sqrt{\xi} q + \epsilon^2 > r^2 .
\]
Therefore
\[
r^2 + r < \xi q^2 + 2  \epsilon \sqrt{\xi} q + \epsilon^2 + 
\frac{\eta + \sqrt{\xi}}{2} q
< \xi q^2 + \eta q + \epsilon^2
< \xi q^2 + \eta q + 2 
\]
by the equation (\ref{**}).

Next we shall prove the assertion (2).
Suppose $\sqrt{\xi} \in {\Bbb Z}$.
Since $r > \sqrt{\xi} q$, we may assume that $r = \sqrt{\xi} q + 1$.
Then,
\[
(\xi q^2 + \eta q + 2) - (r^2 + r)
= (\eta - 3 \sqrt{\xi} ) q .
\]
Therefore, the assertion (2) immediately follows from this.

Now, we shall prove the assertion (3).
Assume $\eta > \sqrt{\xi}$.
Since $(a,b,c) \neq (1,1,1)$, we know $\xi > 1$.
If $\sqrt{\xi} \not\in {\Bbb Z}$, then the assertion immediately 
follows from the assertion (1).
Therefore, we may assume $\sqrt{\xi} \in {\Bbb Z}$.

Let $n$, $q$ and $v$ be integers such that
\[
n = \xi q + v, \ \ v = \sqrt{\xi} -1 .
\]
We set 
\[
r = \sqrt{\xi} q + 1 .
\]
Then,
\[
\sqrt{\xi} r = \xi q + \sqrt{\xi} > n .
\]
Furthermore, by Lemma~\ref{dim}, 
\begin{eqnarray*}
2 \dim_k S_n - (r^2+r) & = & 
\left( \xi q^2 + (\eta + 2v) q + 2 \dim_k S_v \right)
- \left( \xi q^2 + 3 \sqrt{\xi} q + 2 \right) \\
& = & 
(\eta - \sqrt{\xi} -2)q + (2 \dim_k S_v - 2) .
\end{eqnarray*}
Since $\eta - \sqrt{\xi}$ is a non-negative integer, we know
$\eta - \sqrt{\xi} \geq 3$ by Lemma~\ref{012}.
Consequently, we have $2 \dim_k S_n - (r^2+r) > 0$ for $q \gg 0$.
\qed

\begin{Remark}
\begin{rm}
If $(a+b+c)^2 > abc$, then $R_s({\frak p})$ is Noetherian
by a result of Cutkosky~\cite{Cu}.

If $(a+b+c)^2 > abc$ and $\sqrt{abc} \not\in {\Bbb Q}$,
then the existence of a negative curve follows from 
Nakai's criterion for ampleness, Kleimann's theorem and
the cone theorem
(e.g.\ Theorem~1.2.23 and Theorem~1.4.23 in \cite{La}, 
Theorem~4-2-1 in \cite{KMM}).

The condition $(a+b+c)^2 > abc$ is equivalent to 
$(-K_X)^2 > 0$.
If $-K_X$ is ample, then the finite generation of
the total coordinate ring follows from 
Proposition~2.9 and Corollary~2.16 in Hu-Keel~\cite{kee-git}.

If $(a,b,c)=(5,6,7)$, then the negative curve $C$ is the
proper transform of the curve defined by $y^2-zx$.
Therefore, $C$ is linearly equivalent to $12A - E$.
Since $(a+b+c)^2 > abc$, $(-K_X)^2 > 0$.
Since
\[
-K_X.C = (18A-E).(12A-E) = 0.028\cdots > 0,
\]
$-K_X$ is ample by Nakai's criterion.

If $(a,b,c)=(7,8,9)$, then the negative curve $C$ is the
proper transform of the curve defined by $y^2-zx$.
Therefore, $C$ is linearly equivalent to $16A - E$.
Since $(a+b+c)^2 > abc$, $(-K_X)^2 > 0$.
Since
\[
-K_X.C = (24A-E).(16A-E) = -0.23\cdots < 0,
\]
$-K_X$ is not ample by Nakai's criterion.
\end{rm}
\end{Remark}

\section{Degree of a negative curve}\label{sec5}

\begin{Remark}\label{4.1}
\begin{rm}
Let $k$ be a field of characteristic zero, and $R$
be a polynomial ring over $k$ with variables 
$x_1$, $x_2$, \ldots, $x_m$.
Suppose that $P$ is a prime ideal of $R$.
By \cite{LST}, we have
\[
P^{(r)}
= 
\left\{
 h \in R 
\ \left| \
0 \leq \alpha_1 + \cdots + \alpha_m < r 
   \Longrightarrow 
   \frac{
         \partial^{\alpha_1 + \cdots + \alpha_m} h
        }{
          \partial x_1^{\alpha_1} \cdots \partial x_m^{\alpha_m} 
        } \in P
\right.
\right\} .
\]
In particular, if $f \in P^{(r)}$, then
\[
\frac{\partial f}{\partial x_1}, \ldots, 
\frac{\partial f}{\partial x_m} \in P^{(r-1)} .
\]
\end{rm}
\end{Remark}

\begin{Proposition}\label{4.2}
Let $a$, $b$ and $c$ be pairwise coprime integers,
and $k$ be a field of characteristic zero.
Suppose that a negative curve exists, i.e., 
there exist positive integers $n$ and $r$ 
satisfying $[{{\frak p}_k(a,b,c)}^{(r)}]_n \neq 0$ and
$n/r < \sqrt{abc}$.

Set $n_0$ and $r_0$ to be 
\begin{eqnarray*}
n_0 & = &
\min \{ n \in {\Bbb N} \mid 
\mbox{$\exists r>0$ such that $n/r < \sqrt{\xi}$
and $[{\frak p}^{(r)}]_n \neq 0$} \} \\
r_0 & = & \lfloor \frac{n_0}{\sqrt{\xi}} \rfloor + 1 ,
\end{eqnarray*}
where $\lfloor \frac{n_0}{\sqrt{\xi}} \rfloor$ is the
maximum integer which is less than or equal to $\frac{n_0}{\sqrt{\xi}}$.

Then, the negative curve $C$ is linearly equivalent to
$n_0A - r_0E$.
\end{Proposition}

\proof
Suppose that the negative curve $C$ is linearly equivalent to
$n_1A - r_1E$.
Since $n_1/r_1 < \sqrt{\xi}$ and 
$[{\frak p}^{(r_1)}]_{n_1} \neq 0$, we have $n_1 \geq n_0$.
Since $H^0(X, {\cal O}(n_0A - r_0E)) \neq 0$ with 
$n_0/r_0 < \sqrt{abc}$,
$n_0A - r_0E -C$ is linearly equivalent to an effective divisor.
Therefore, $n_0 \geq n_1$.
Hence, $n_0 = n_1$.

Since $n_0/r_1 < \sqrt{\xi}$, 
$r_0 \leq r_1$ holds.
Now, suppose $r_0 < r_1$.
Let $f$ be the defining equation of $\pi(C)$,
where $\pi : X \rightarrow {\Bbb P}$ is the blow-up
at $V_+({\frak p})$.
Then, we have
\[
[{\frak p}^{(r_1-1)}]_{n_0} = [{\frak p}^{(r_1)}]_{n_0}
= k \ f .
\]
If $n$ is an integer less than $n_0$, then
$[{\frak p}^{(r_1-1)}]_n = 0$ because
\[
\frac{n}{r_1-1} < \frac{n_0}{r_1-1} \leq \frac{n_0}{r_0} < \sqrt{\xi} .
\]
%
By Remark~\ref{4.1}, we have
\[
\frac{\partial f}{\partial x}, \frac{\partial f}{\partial y}, 
\frac{\partial f}{\partial z} \in {\frak p}^{(r_1-1)} .
\]
Since their degrees are strictly less than $n_0$,
we know
\[
\frac{\partial f}{\partial x} = \frac{\partial f}{\partial y} = 
\frac{\partial f}{\partial z} = 0 .
\]
On the other hand, the equality
\[
ax \frac{\partial f}{\partial x} +
by \frac{\partial f}{\partial y} +
cz \frac{\partial f}{\partial z} 
= n_0 f 
\]
holds.
Remember that $k$ is a field of
characteristic zero.
It is a contradiction.
\qed

\begin{Remark}\label{rem5.3}
\begin{rm}
Let $a$, $b$ and $c$ be pairwise coprime integers,
and $k$ be a field of characteristic zero.
Assume that the negative curve $C$ exists, and
$C$ is linearly equivalent to
$n_0A - r_0E$.

Then, by Proposition~\ref{4.2},
we obtain
\begin{eqnarray*}
n_0 & = &
\min \{ n \in {\Bbb N} \mid 
[{\frak p}^{(\lfloor \frac{n}{\sqrt{\xi}} \rfloor + 1)}]_n \neq 0  \} \\
r_0 & = & \lfloor \frac{n_0}{\sqrt{\xi}} \rfloor + 1 .
\end{eqnarray*}
\end{rm}
\end{Remark}

\begin{Theorem}\label{bound}
Let $a$, $b$ and $c$ be pairwise coprime positive integers such that
$\sqrt{\xi} > \eta$.
Assume that (C2) is satisfied, i.e., there exist positive integers
$n_1$ and $r_1$ such that $n_1/r_1 < \sqrt{\xi}$ and
$\dim_k S_{n_1} > r_1(r_1+1)/2$.
Suppose $n_1 = \xi q_1 + v_1$, where $q_1$ and $v_1$ are integers 
such that $0 \leq v_1 < \xi$.

Then, $q_1 < \frac{2 \dim_k S_{v_1}}{\sqrt{\xi} - \eta}$ holds.

In particular,
\[
n_1 = \xi q_1 + v_1
< \frac{2 \xi \max\{ \dim_k S_t \mid 0 \leq t < \xi \}}{\sqrt{\xi} - \eta}
+ \xi .
\]
\end{Theorem}

\proof
We have
\[
r_1 > \frac{n_1}{\sqrt{\xi}} = \sqrt{\xi} q_1 + \frac{v_1}{\sqrt{\xi}} .
\]
Therefore,
\[
2 \dim_k S_{n_1} > 
r_1^2 + r_1 >
\xi q_1^2 + 2v_1q_1 + \frac{v_1^2}{\xi} + \sqrt{\xi} q_1 + 
\frac{v_1}{\sqrt{\xi}} .
\]
By Lemma~\ref{dim}, we have
\[
(\sqrt{\xi} - \eta) q_1 < 2 \dim_k S_{v_1} - \frac{v_1^2}{\xi} - 
\frac{v_1}{\sqrt{\xi}} \leq 2 \dim_k S_{v_1} .
\]
\qed

Remember that, if $\sqrt{\xi} < \eta$, then (C2) is always satisfied by Theorem 4.3 (3).

\section{Calculation by computer}\label{sec6}

In this section, we assume that the characteristic of $k$ is zero.

\subsection{Examples that do not satisfy (C2)}\label{c2}

Suppose that (C2) is satisfied, i.e., there exist positive integers
$n_1$ and $r_1$ such that $n_1/r_1 < \sqrt{\xi}$ and
$\dim_k S_{n_1} > r_1(r_1+1)/2$.
Put $n_1 = \xi q_1 + v_1$, where $q_1$ and $v_1$ are integers 
such that $0 \leq v_1 < \xi$.
If $\sqrt{\xi} > \eta$, then 
$q_1 < \frac{2 \dim_k S_{v_1}}{\sqrt{\xi} - \eta}$ holds
 by Theorem~\ref{bound}.

By the following programming on MATHEMATICA,
we can check whether (C2) is satisfied or not
in the case where $\sqrt{\xi} > \eta$.

\begin{verbatim}
c2[a_, b_, c_] := 
Do[
  If[(a + b + c)^2 > a b c , Print["-K: self-int positive"]; Goto[fin]];
  s = Series[((1 - t^a)(1 - t^b)(1 - t^c))^(-1), {t, 0, a b c}];
  Do[   h = SeriesCoefficient[s, k];
        m = IntegerPart[2 h/(Sqrt[a b c] - a - b - c)];
        Do[  r = IntegerPart[(a b c q + k)(Sqrt[a b c]^(-1))] + 1;
             If[2 h + q(a + b + c) + a b c q^2 + 2q k > r (r + 1), 
               Print[StringForm["n=``, r=``", a b c q + k, r]];
               Goto[fin]], 
          {q, 0, m}], 
    {k, 0, a b c - 1}];
  Print["c2 is not satisfied"];     
  Label[fin];
  Print["finished"]]
\end{verbatim}

Calculations by a computer show that (C2) is not satisfied in some cases,
for example, $(a,b,c) = (5,33,49), (7,11,20), (9,10,13), \cdots$.

The examples due to Goto-Nishida-Watanabe~\cite{GNW}
 have negative curves
with $r = 1$.
Therefore, by Remark~\ref{r=01}, they satisfy the condition (C2).

In the case where $(a,b,c) = (5,33,49), (7,11,20), (9,10,13), \cdots$, 
the authors do not know whether $R_s({\frak p}_k)$ is Noetherian or not.

\begin{Remark}
\begin{rm}
Set
\begin{eqnarray*}
A&=& \{(a,b,c) \mid 0<a\le b\le c \le 50, a,b,c \mbox{ are pairwise coprime}\}\\
B&=& \{(a,b,c) \in A \mid a+b+c > \sqrt{abc}\}\\
C&=& \{(a,b,c) \in A \mid (a,b,c) \mbox{ does not satisfy (C2)}\}.
\end{eqnarray*}
$\sharp A = 6156$, $\sharp B = 1950$, $\sharp C = 457$.
By Theorem 4.3, we know $B \cap C = \emptyset$.
\end{rm}
\end{Remark}

\subsection{Does a negative curve exist?}

By the following simple computer programming on MATHEMATICA, it is possible to know whether a negative curve exists or not.

\begin{verbatim}
n[a1_, b1_, c1_, r1_, d1_] := (V = 0;
  Do[
   mono = {};
   Do[ e1 = d1 - i*a1; 
     Do[ h1 = e1 - j*b1; k1 = Floor[h1/c1]; 
       If[ h1 / c1 == k1, 
         mono = Join[mono, {x^i y^j z^(k1)}]], 
         {j, 0, Floor[e1/b1]}
       ], {i, 0, Floor[d1/a1]}
     ];
   w = Length[mono];
   If[w > N[r1*(r1 + 1)/2], 
     V = 1; W = w; J = d1; R = r1; H = N[r1*(r1 + 1)/2], 
     If[ w > 0,
       f[x_, y_, z_] := mono;
       mat = {};
       Do[
         Do[
           mat = Join[mat, { D[f[x, y, z], {x, j}, {y, i - j}] }], {j, 0, i}
           ], {i, 0, r1 - 1}
         ];
       mat = mat /. x -> 1 /. y -> 1 /. z -> 1;
       q = MatrixRank[mat]; 
       If[ q < w, V = 1; W = w; J = d1; R = r1; H = N[r1*(r1 + 1)/2] ]
      ]
    ]
  ]);
\end{verbatim}

\begin{verbatim}
t[a_, b_, c_, rmade_] := (
  Do[
    W = 0;
    p = Ceiling[r*Sqrt[a*b*c]] - 1;
    Do[
      n[a, b, c, r, p - u]; 
      If[V == 1, 
        J1 = J; Break[]], {u, 0, a - 1}
      ];
    If[
      V == 1, 
      Do[
        n[a, b, c, r, J1 - a*u]; 
        If[V == 0, 
          J1 = J; Break[]], {u, 1, b*c}
        ];

      Do[
        n[a, b, c, r, J1 - b*u]; 
        If[V == 0, 
          J1 = J; Break[]], {u, 1, a*c}
        ];

      Do[
        n[a, b, c, r, J1 - c*u]; 
        If[V == 0, 
          J1 = J; Break[]], {u, 1, c*a}
        ]
      ];
    If[W > 0, Break[]];

    Print["r th symbolic power does not contain a negative curve if r <= ",
      r], {r, 1, rmade}
    ];
  If[W == 0, 
    Print["finished"], 
    Print["There exists a negative curve. Degree = ", J, ",  r = ", R, 
      ",  Dimension of homog. comp. = ", W, ",  # of equations = ", H]
    ]
  )
\end{verbatim}

By the command t[a,b,c,r], we can check whether ${\frak{p}}(a,b,c)^{(m)}$ contains an equation of a negative curve for $m = 1,2, \cdots , r$.


${\frak{p}}(9,10,13)^{(m)}$ does not contain an equation of a negative curve if $m \le 24$. 
Remember that $(9,10,13)$ does not satisfy (C2).
Our computer gave up computation of ${\frak{p}}(9,10,13)^{(25)}$ for scarcity of memories. We don't know whether ther exists a negative curve in the case $(9,10,13)$.

On the other hand, there are examples that (C2) is not satisfied but there exists a negative curve.

\begin{itemize}
\item Suppose $(a,b,c) = (5,33,49)$. Then (C2) is not satisfied, but $[{\frak{p}}(5,33,49)^{(18)}]_{1617}$ contatins a negative curve.
\item Suppose $(a,b,c) = (8,15,43)$. Then (C2) is not satisfied, but 
$[{\frak{p}}(8,15,43)^{(9)}]_{645}$ contains a negative curve.
\end{itemize}

\noindent
\begin{tabular}{l}
Kazuhiko Kurano \\
Department of Mathematics \\
Faculty of Science and Technology \\
Meiji University \\
Higashimita 1-1-1, Tama-ku \\
Kawasaki 214-8571, Japan \\
{\tt kurano@math.meiji.ac.jp} \\
{\tt http://www.math.meiji.ac.jp/\~{}kurano}
\end{tabular}

\vspace{2mm}

\noindent
\begin{tabular}{l}
Naoyuki Matsuoka \\
Department of Mathematics \\
Faculty of Science and Technology \\
Meiji University \\
Higashimita 1-1-1, Tama-ku \\
Kawasaki 214-8571, Japan \\
{\tt matsuoka@math.meiji.ac.jp} 
\end{tabular}

\end{document}